\documentclass[12pt]{article}

\usepackage{amsmath,amsthm,amsfonts,amssymb,amscd}

\theoremstyle{plain}
\newtheorem{thm}{Theorem}[section]
\newtheorem{rk}[thm]{Remark}
\newtheorem{prop}[thm]{Proposition}
\newtheorem{clly}[thm]{Corollary}
\newtheorem{lemma}[thm]{Lemma}

\newtheorem{maintheorem}{Theorem}
\newtheorem{claim}[thm]{Claim}

\newcommand{\pf}{{\flushleft{\bf Proof: }}}

\newcommand{\re}{\mathbb{R}}

\newcommand{\rt}{\rightarrow}

\newcommand{\dist}{\operatorname{dist}}

\newcommand{\interior}{\operatorname{int}}

\newcommand{\cl}{\operatorname{Cl}} 

\newcommand{\di}{\operatorname{dim}}

\newcommand{\per}{\operatorname{Per}}

\newcommand{\cri}{\operatorname{Crit}}

\title{Homoclinic classes for generic $C^1$ vector
fields}

\author{C. M. Carballo, C. A. Morales \& M. J. Pacifico\thanks
{This work is partially supported by CNPq, FAPERJ and PRONEX/Dyn. Systems}}

\date{June 1, 2000\\
      Revised July 5, 2001}

\begin{document}

\maketitle

\begin{abstract}
We prove that
homoclinic classes for a {\em  residual} set of $C^1$ vector fields
$X$ on closed $n$-manifolds are maximal transitive,
and depend continuously on periodic orbit data.
In addition, $X$ does not exhibit
cycles formed by homoclinic classes.
We also prove that a homoclinic class 
of $X$ is isolated if and only if it is
$\Omega$-isolated, and it is the intersection of
its stable set with its unstable set.
All these properties are well known for
structural stable Axiom A vector fields.
\end{abstract}

\section{Introduction}
\label{s-intro}

We show some properties of homoclinic classes for generic
$C^1$ flows on closed $n$-manifolds.
By homoclinic class we mean the closure
of the transversal
homoclinic points associated to a
hyperbolic periodic orbit.
So, homoclinic classes are transitive and
the closure of its periodic orbits \cite[Chapter 2, $\S 8$]{PT93}.

For structural stable Axiom A vector fields it is
known that homoclinic classes are maximal transitive
and depend continuously on the periodic orbit data.
In addition, if $H$ is a homoclinic class of $X$ then it is saturated,
that is, $H=W_X^s(H) \cap W_X^u(H)$, where $W_X^s(H)$ is the stable set of
$H$ and $W^u_X(H)$ is the unstable set of $H$ \cite[p. 371]{BBV99}.
Moreover, such vector fields do not exhibit
cycles formed by homoclinic classes.
In this paper we shall prove these
properties for generic $C^1$ vector fields
on closed $n$-manifolds $M$ (neither structural stability nor Axiom A 
is assumed).
Furthermore, we prove that generically a
homoclinic class is isolated if and only if it is
isolated from the nonwandering
set.
In particular, all the mentioned properties hold
for a dense set of $C^1$ vector fields on $M$.
It is interesting to observe that neither structural stability nor Axiom A
is dense in the space of $C^1$ vector fields on $M$, $\forall n\geq 3$.

To state our results in a precise way we use the following notation.
$M$ is a compact boundaryless $n$-manifold, ${\cal X}^1(M)$ is the
space of $C^1$ vector fields endowed with the $C^1$ topology.
Given $X\in {\cal X}^1(M)$, $X_t$ denotes the flow induced by $X$.
The $\omega$-limit set of $p$ 
is the set $\omega_X(p)$ of accumulation points of the 
positive orbit of $p$.
The $\alpha$-limit set of $p$ is
$\alpha_X(p)=\omega_{-X}(p)$,
where $-X$ denotes the time-reversed flow
of $X$.
The nonwandering set $\Omega(X)$ of $X$
is the set of $p$ such that for every neighborhood
$U$ of $p$ and $T>0$ there is $t>T$ such that
$X_t(U)\cap U\neq\emptyset$. Clearly $\Omega(X)$ is
closed, nonempty and contains any $\omega$-limit
($\alpha$-limit) set. 
A compact invariant set $B$ of $X$ is
{\em $\Omega$-isolated} if
$\Omega(X)\setminus B$ is closed. $B$ is
{\em isolated}
if $B=\cap_{t\in \re} X_t(U)$
for some compact neighborhood $U$ of $B$
(in this case $U$ is called {\em isolating block}).
We denote by $\per(X)$ the union of the periodic
orbits of $X$ and $\cri(X)$ the set formed by the union of $\per(X)$
and the singularities of $X$.

A set is {\em transitive} for $X$ if
it is the $\omega$-limit set of one of its orbits.
A transitive set $\Lambda$ of $X$ is {\em maximal transitive} 
if it contains every transitive
set $T$ of $X$ satisfying $\Lambda\cap T\neq\emptyset$.
Note that a maximal transitive set is maximal with respect to
the inclusion order.
In \cite{BD99, BDP99} it was asked whether every 
homoclinic class $H_f(p)$ of a generic diffeomorphism
$f$ satisfies the property that if $T$ is a transitive set
of $f$ and $p\in T$, then $T\subset H_f(p)$. 
In \cite{Ar99}, M. C. Arnaud also considered homoclinic classes
for $C^1$ diffeomorphisms on $M$,
and in particular she gives a positive answer to this question 
\cite[Corollary 40]{Ar99}.
On the other hand, item (1) of Theorem \ref{thA} below states 
that generically any transitive set of a $C^1$ vector field 
intersecting the homoclinic class is included in it, and thus
the diffeomorphism version of it extends this result of M. C. Arnaud.

If $\Lambda$ is a compact invariant set of $X$, we denote
\[
W^s_X(\Lambda)=\{q\in M: \dist(X_t(q),\Lambda)\to 0,
t\to\infty \}
\]
and
\[
W^u_X(\Lambda)=\{q\in M: \dist(X_t(q),\Lambda)
\to 0, t\to -\infty \},
\]
where $\dist$ is the metric on $M$.
These sets are called respectively the stable and unstable 
set of $\Lambda$.
We shall denote $W^s_X(p)=W^s_X({\cal O}_X(p))$
and $W^u_X(p)=W^u_X({\cal O}_X(p))$
where ${\cal O}_X(p)$ is the orbit of $p$.
We say that $\Lambda$ is {\em saturated} if 
$W^s_X(\Lambda)\cap W^u_X(\Lambda)= \Lambda$.

A {\em cycle} of $X$ is a finite set of 
compact invariant sets $\Lambda_0,\Lambda_1,\dots,\Lambda_n$ 
such that $\Lambda_{n}=\Lambda_0$, and
$\Lambda_0,\Lambda_1,\dots,\Lambda_{n-1}$ are disjoint, and
$$
(W^u_X(\Lambda_i)\setminus \Lambda_{i}) \cap 
(W^s_X(\Lambda_{i+1})\setminus \Lambda_{i+1}) \neq \emptyset
$$
for all $i=0,...,n-1$.

A compact invariant set $\Lambda$
of $X$ is {\em hyperbolic}
if there is a continuous tangent
bundle decomposition $E^s\oplus E^X\oplus E^u$
over $\Lambda$ such that $E^s$ is contracting,
$E^u$ is expanding and $E^X$ denotes
the direction of $X$.
We say that $p\in \cri(X)$ is hyperbolic
if ${\cal O}_X(p)$
is a hyperbolic set
of $X$.

The Stable Manifold Theorem \cite{HPS70}
asserts that $W^s_X(p)$ is an immersed manifold
tangent to $E^s\oplus E^X$ for every $p$
in a hyperbolic set $\Lambda$ of $X$. Similarly for $W^u_X(p)$. 
This remark applies when $\Lambda={\cal O}_X(p)$ for some
$p\in \cri(X)$ hyperbolic.
As already defined, the homoclinic class associated
to a hyperbolic periodic orbit $p$ of $X$, $H_X(p)$,
is the closure of the transversal intersection
orbits in $W^s_X(p)\cap W^u_X(p)$.

We say that $X$ is {\em Axiom A} if
$\Omega(X)$ is both hyperbolic and
the closure of $\cri(X)$.
The non wandering set
of a nonsingular Axiom A flow splits in a finite disjoint
union of homoclinic classes \cite[Chapter 0, p. 3]{PT93}.

Another interesting property of homoclinic classes
for Axiom A vector fields is their continuous dependence
on the periodic orbit data,
that is, the map $p\in\per (X) \to H_X(p)$
is upper-semicontinous.

In general, we say that a compact set sequence
$\Lambda_n$ {\em accumulates} on a compact set
$\Lambda$ if for every neighborhood $U$ of $\Lambda$
there is $n_0>0$ such that $\Lambda_n\subset U$
for all $n\geq n_0$.
Note that this kind of accumulation is weaker
than the usual Hausdorff metric accumulation.

If ${\cal Y}$ denotes a metric space,
then ${\cal R}\subset {\cal Y}$ is {\em residual} in ${\cal Y}$ if
${\cal R}$ contains a countable intersection of open-dense subsets 
of ${\cal Y}$.
Clearly a countable intersection of residual subsets of ${\cal Y}$
is a residual subset of ${\cal Y}$.
For example, the set of Kupka--Smale vector fields ${\cal KS}^1(M)$ 
on $M$ is a residual subset of ${\cal X}^1(M)$ 
\cite[Chapter 3, $\S$ 3]{MP92}.
Recall that a vector field is Kupka--Smale if
all its periodic orbits and singularities are hyperbolic and
the invariant manifolds of such elements intersect transversally.

\begin{maintheorem}
\label{thA}
The following properties hold
for a residual subset of vector fields $X$ in ${\cal X}^1(M)$:
\begin{enumerate}
\item
The homoclinic classes of $X$ are
maximal transitive sets of $X$.
In particular,
different homoclinic classes of $X$ are disjoint.
\item
The homoclinic classes of $X$ are saturated.
\item
The homoclinic classes of $X$ depends
continuously on the periodic orbit data,
that is, the map $p\in\per (X) \to H_X(p)$
is upper-semicontinous.
\item
A homoclinic class of $X$ is isolated if and only if it is $\Omega$-isolated.
\item
The hyperbolic homoclinic classes of $X$ are
isolated.
\item
There is no cycle of $X$ formed by
homoclinic classes of $X$.
\item $X$ has finitely many homoclinic
classes if and only if
the union of the homoclinic classes of $X$
is closed
and every homoclinic class of $X$ is
isolated. 
\end{enumerate}
\end{maintheorem}

When $M$ has dimension three we obtain
the following corollaries using
Theorem \ref{thA}, \cite{Li83}, and \cite{Man82}.
Recall that an isolated set $\Lambda$
of a $C^r$ vector field $X$
is $C^r$ robust transitive ($r\geq 1$) if
it exhibits an isolating block
$U$ such that, for every vector field
$Y$ $C^r$ close to $X$,
$\cap_{t\in I\!\! R}Y_t(U)$
is both transitive and nontrivial
for $Y$.

\begin{clly}
\label{c2}
The properties below are equivalent
for a residual set of nonsingular 3-dimensional
$C^1$ vector fields $X$ and every nontrivial
homoclinic class $H_X(p)$ of $X$:

\begin{enumerate}
\item
$H_X(p)$ is hyperbolic.
\item 
$H_X(p)$ is
isolated.
\item
$H_X(p)$ is
$C^1$ robust transitive for $X$.
\end{enumerate}
\end{clly}

\begin{clly}
\label{c5}
The properties below
are equivalent for a residual set of
nonsingular 3-dimensional $C^1$ vector fields $X$:
\begin{enumerate}
\item
$X$ is Axiom A.
\item
$X$ has finitely many homoclinic classes.
\item
The union of the homoclinic classes
of $X$ is closed
and every homoclinic class of $X$ is isolated.
\end{enumerate}
\end{clly}

The equivalence between the Items (1) and (2) of the above
corollary follows from \cite{Li83, Man82}. It shows
how difficult is to prove the genericity
of vector fields exhibiting finitely many homoclinic classes.
The equivalence between (2) and (3) follows
from Theorem \ref{thA}-(7).

To prove Theorem \ref{thA}
we show in Section \ref{s-hc} that homoclinic
classes $H_X(p)$ for a residual set of $C^1$ vector fields
$X$ satisfy $H_X(p)=\Lambda^+\cap\Lambda^-$,
where $\Lambda^+$ is Lyapunov stable
for $X$ and $\Lambda^-$ is Lyapunov stable for $-X$.
The main technical tool to prove such result is 
Lemma \ref{lwx}, a stronger version of Hayashi's
$C^1$ Connecting Lemma \cite{Ha97}, recently published in 
\cite[Theorem E, p. 5214]{WX00}
(see also \cite{Ar99, Ha1, Ha98}).
In Section \ref{s-ls} we study compact invariant sets
$\Lambda$ of $X$
satisfying $\Lambda=\Lambda^+\cap \Lambda^-$, where
$\Lambda^\pm$ is Lyapunov stable for $\pm X$.
The proof of Theorem \ref{thA} and Corollary \ref{c2}
will be given in the final section
using the results of Sections \ref{s-ls} and \ref{s-hc}.

\begin{rk}
We observe that Theorem \ref{thA} is valid for a residual set of
$C^1$ diffeomorphisms on any $n$-manifold $M$ by the usual method 
of suspension.
\end{rk}

We are thankful to S. Hayashi for useful conversations.

\section{Lyapunov stability lemmas}
\label{s-ls}

In this section we shall establish useful
properties of Lyapunov stable sets.
A reference for Lyapunov stability theory
is \cite{BS70}.

Recall we have denoted by $X_t$, $t\in I\!\! R$ 
the flow generated by $X\in {\cal X}^1(M)$.
Given $A\subset M$ and $R\subset I\!\! R$
we set
$X_R(A)=\{X_t(q):(q,t)\in A\times R\}$.
We denote $\cl(A)$ the closure of $A$, and $\interior(A)$
the interior of $A$.
If $\epsilon>0$ and
$q\in M$ we set
$B_\epsilon(q)$ the $\epsilon$-ball centered at $q$.

A compact subset $A \subseteq M$ is {\em Lyapunov stable} for 
$X$ if for every open set $U$ containing $A$
there exists an open set $V$ containing $A$
such that $X_t(V)\subset U$ for every $t \geq 0$.
Clearly a Lyapunov stable set
is forward invariant.

The following lemma summarizes
some classical properties of Lyapunov stable sets
(see \cite[Chapter V]{BS70}).

\begin{lemma}
\label{l0}
Let $\Lambda^+$ be a Lyapunov stable set of $X$.
Then,
\begin{enumerate}
\item
If $x_n\in M$ and $t_n\geq 0$ satisfy
$x_n \rt x\in \Lambda^+$ and $X_{t_n}(x_n) \rt y$,
then
$y\in \Lambda^+$;
\item
$W^u_X(\Lambda^+)\subset \Lambda^+$;
\item
if $\Gamma$ is a transitive set of $X$
and $\Gamma\cap\Lambda^+\neq\emptyset$, then
$\Gamma\subset \Lambda^+$.
\end{enumerate}
\end{lemma}

We are interested in invariant compact sets
$\Lambda=\Lambda^+\cap\Lambda^-$ of $X$,
where
$\Lambda^+$ is Lyapunov stable set for $X$ and 
$\Lambda^-$ is Lyapunov stable set for the reversed flow $-X$.
We shall call such sets {\em neutral} for the sake of simplicity.
As we shall see in the next section, homoclinic classes
are neutral sets for generic $C^1$ vector fields on closed $n$-manifolds.

Elementary properties of neutral sets are given
in the lemma below.

\begin{lemma} 
\label{l1}
Let $\Lambda$ be a neutral set
of $X$. Then,
\begin{enumerate}
\item
$\Lambda$ is saturated;
\item
$\Lambda$ is transitive for $X$ if and only if
$\Lambda$ is maximal transitive
for $X$. In particular, different
transitive neutral sets of $X$ are disjoint.
\end{enumerate}
\end{lemma}

{\pf}
Let $\Lambda=\Lambda^+\cap\Lambda^-$
with $\Lambda^\pm$ being Lyapunov stable for $\pm X$.
Clearly $W^u_X(\Lambda)\subset
\Lambda^+$ by Lemma \ref{l0}-(2).
Similarly, $W^s_X(\Lambda)\subset
\Lambda^-$.  
Hence
$$
W^u_X(\Lambda)\cap W^s_X(\Lambda)\subset
\Lambda^+\cap \Lambda^-=\Lambda.
$$
Conversely, $\Lambda\subset
W^u_X(\Lambda)\cap W^s_X(\Lambda)$
since $\Lambda$ is invariant. This proves (1).

Now, by Lemma \ref{l0}-(3),
if $\Gamma$ is a transitive set intersecting
$\Lambda$, then
$\Gamma\subset \Lambda^+$ and $\Gamma\subset
\Lambda^-$.
Thus, $\Gamma\subset\Lambda^+\cap\Lambda^-=\Lambda$,
and so, $\Lambda$ is maximal
transitive. The converse is obvious.
Different transitive neutral sets
of $X$ are maximal transitive, and so,
they are necessarily disjoint. This finishes the proof.
\qed

Note that a Smale horseshoe
with a first tangency is an example of a maximal
transitive set which is not neutral, see Proposition \ref{p2}.
This example also provides a
hyperbolic homoclinic class which is not neutral
(compare with Theorem \ref{h-neutral}).

\begin{prop}
\label{p3}
There is no cycle of $X$ formed by transitive neutral sets. 
\end{prop}

{\pf}
By contradiction suppose that there exists a cycle
$\Lambda_0,\dots , \Lambda_n$ of $X$
such that every $\Lambda_i$ is a transitive
neutral set of $X$. Recall $\Lambda_n=\Lambda_0$.

Set 
$\Lambda_i = \Lambda_i^+ \cap \Lambda_i^-$
where each $\Lambda_i^\pm$ is Lyapunov stable
for $\pm X$.
Choose
$$
x_i \in (W^u_X(\Lambda_i)\setminus
\Lambda_i) \cap (W^s_X(\Lambda_{i+1})
\setminus \Lambda_{i+1})
$$
according to the definition.

We claim that $x_i\in \Lambda_0^-$ for every $i$.
Indeed, as $W^s_X(\Lambda_0)\subset \Lambda_0^-$
one has $x_{n-1} \in \Lambda_0^-$.
Assume by induction
that $x_i\in \Lambda_0^-$
for some $i$.
As
$x_{i} \in W^u_X(\Lambda_{i})$,
the backward invariance of $\Lambda_0^-$
implies
$$
\Lambda_0^-
\cap \Lambda_{i}\supset \alpha_X(x_{i})
\neq\emptyset.
$$
By Lemma \ref{l0}-(3) one has
$\Lambda_{i}\subset \Lambda_0^-$
since $\Lambda_{i}$ is transitive.
In particular, $W^s_X(\Lambda_{i})
\subset \Lambda_0^-$ 
by Lemma \ref{l0}-(2)
applied to $-X$. As $x_{i-1}\in
W^s_X(\Lambda_{i})$, one has
$x_{i-1}\in \Lambda_0^-$.
The claim follows by induction.

By the claim $x_0 \in \Lambda_0^-$.
As $W^u_X(\Lambda_0)\subset \Lambda_0^+$
and $x_0\in W^u_X(\Lambda_0)$
(by definition)
one has $x_0\in\Lambda_0^+\cap \Lambda_0^-=
\Lambda_0$. This contradicts
$x_0\in W^u_X(\Lambda_0)\setminus \Lambda_0$
and the proposition is proved.
\qed

\begin{lemma}
\label{l2}
If $\Lambda$ is neutral for $X$, then
for every neighborhood $U$ of $\Lambda$
there exists a neighborhood $V\subset U$
of $\Lambda$ such that
$$
\Omega(X)\cap V\subset \cap_{t\in I\!\! R}X_t(U).
$$
\end{lemma}

{\pf}
Let $U$ be a neighborhood of
a neutral set $\Lambda$ of $X$.
Choose $U'\subset \cl(U')\subset U$
with $U'$ being another neighborhood of
$\Lambda$.
We claim that there is a neighborhood
$V\subset U'$ of $\Lambda$ so that:

\begin{description}
\item{(1)}
$t\geq 0 \text{ and } p\in V\cap X_{-t}(V) \implies X_{[0,t ]}(p)\subseteq U'$.
\item{(2)}
$t\leq 0 \text{ and } p\in V\cap X_{t}(V) \implies X_{[-t,0]}(p)\subseteq U'$.
\end{description}

Indeed, it were not
true then there would exist a neighborhood $U$
of $\Lambda$ and sequences
$p_n\rt \Lambda$,
$t_n> 0$ such that
$X_{t_n}(p_n)\rt \Lambda$
but
$X_{[0,t_n]}(p_n)\not\subseteq U'$.
Choose $q_n\in X_{[0,t_n]}(p_n)\setminus U'$.
Write $q_n= X_{t_n'}(p_n)$ for some $t_n'\in [0,t_n]$ and
assume that $q_n\rt q$ for some $q\notin U'$.
Let $\Lambda= \Lambda^+\cap \Lambda^-$
with $\Lambda^\pm$ Lyapunov stable for $\pm X$.
Since $\Lambda^+$ is Lyapunov stable for $X$
and $t_n'> 0$, Lemma \ref{l0}-(1) implies $q\in \Lambda^+$.
On the other hand, as we can write
$q_n= X_{t_n'-t_n}(X_{t_n}(p_n))$ where
$t_n'-t_n> 0$ and $X_{t_n}(p_n)\rt \Lambda$ and
using again Lemma \ref{l0}-(1) we have that
$q\in \Lambda^-$.
This proves that $q\in \Lambda$,
a contradiction since $q\notin U'$.
This proves the claim.

Next we prove that
$\Omega(X)\cap V\subseteq \cap_{t\in \re}X_t(U)$.
Indeed,
choose $q\in \Omega(X)\cap V$.
By contradiction,
we assume that there is $t_0>0$ (say) such that $X_{t_0}(q)\notin U$.
Then, there is a ball $B_{\epsilon}(q)\subseteq V$ such that
$X_{t_0}(B_{\epsilon}(q))\cap \cl(U')=\emptyset$.
As $q\in \Omega(X)$ there exists
$t> t_0$ such that 
$B_{\epsilon}(q)\cap X_{-t}(B_{\epsilon}(q))\neq \emptyset$.
Pick $p\in B_\epsilon(q)\cap X_{-t}(B_\epsilon(q))$.
By (1) above one has $X_{[0,t]}(p)\subseteq U'$
since $B_\epsilon(q)\subset V$.
This contradicts
$X_{t_0}(p)\in X_{t_0}(B_{\epsilon}(q))$
and $X_{t_0}(B_{\epsilon}(q))\cap \cl(U')=\emptyset$.
The proof is completed.
\qed

A first consequence of the above lemma is
the following corollary.
Given compact subsets $A,B \subset M$ we denote $\dist (A,B)=
\inf \{\dist(a,b); a \in A, b \in B\}$.

\begin{clly}
\label{cor1}
If $\Lambda$ is a neutral set of $X$ and
$\Lambda_n$ is a sequence of transitive sets
of $X$ such that $\dist(\Lambda_n,\Lambda)\to 0$ as
$n\to \infty$, then $\Lambda_n$ accumulates on $\Lambda$.
\end{clly}

{\pf}
Let $\Lambda_n$ and $\Lambda$ as in the statement.
Fix a neighborhood $U$ of $\Lambda$
and let $V\subset U$ be the neighborhood of $\Lambda$ obtained by the previous lemma.
As $\dist(\Lambda_n,\Lambda)\to 0$ as $n\to\infty$
we have that $\Lambda_n\cap V\neq\emptyset$ for every $n$ large.
Let $q_n$ the dense orbit of $\Lambda_n$.
Clearly $q_n\in \Omega(X)$.
We can assume that $q_n\in V$ for $n$ large,
and so, $q_n\in \Omega(X)\cap V$.
Then, $X_t(q_n)\in U$
for every $t$. In particular, $\Lambda_n=\omega_X(q_n)\subset
\cl(U)$. This proves the corollary since
$U$ is arbitrary.
\qed

\begin{prop}
\label{p2}
A neutral set is isolated if and only if
it is $\Omega$-isolated.
\end{prop}

{\pf}
We first claim that any
saturated $\Omega$-isolated set $\Lambda$ of $X$
is isolated.
Indeed, since $\Lambda$ is
$\Omega$-isolated, there is
$U\supset \Lambda$ open such that
$\cl(U)\cap \Omega(X)= \Lambda$.
This $U$ is an isolating block for $\Lambda$.
For if $x\in \cap_{t\in I\!\! R}X_t(U)$, then
$\omega_X(x)\cup \alpha_X(x)\subset \cl(U)
\cap \Omega(X)\subset \Lambda$. So, $x\in
W^s_X(\Lambda)\cap W^u_X(\Lambda)=\Lambda$.
This proves that $\cap_{t\in I\!\! R}X_t(U)
\subset \Lambda$. The opposite inclusion follows
since $\Lambda$ is invariant. The claim follows.

To prove that invariant $\Omega$-isolated
neutral set are isolated we use the above claim
and Lemma \ref{l1}-(1).
To prove that
isolated neutral sets are $\Omega$-isolated
we use Lemma \ref{l2}.
\qed

\begin{prop}
\label{p2.5}
Transitive hyperbolic neutral sets are isolated.
\end{prop}

{\pf}
By Proposition \ref{p2} it is suffices to show that transitive
neutral hyperbolic
sets $\Lambda$ are $\Omega$-isolated.

Suppose by contradiction that
$\Lambda$ is not $\Omega$-isolated.
Then, there is a sequence
$p_n\in \Omega(X)\setminus \Lambda$
converging to $p\in \Lambda$.
Fix $U$ a neighborhood of $\Lambda$ and
let $V$ be given in Lemma \ref{l2} for
$U$. We can assume that $p_n\in V$ for every $n$.
As $p_n$ is non wandering for $X$,
for every $n$ there are sequences
$q_i\in V\to p_n$ and
$t_i>0$ such that $X_{t_i}(q_i)\to p_n$
as $i\to \infty$.
By (1) in the proof of Lemma \ref{l2} we have
$X_{[0,t_i]}(q_i)\subset U$
for every $i$.
So, we can construct a
periodic pseudo orbit of $X$ arbitrarily close to $p_n$.
By the Shadowing Lemma for Flows 
(\cite[Theorem 18.1.6, p. 569]{KH95})
applied to the hyperbolic set $\Lambda$,
such a periodic pseudo orbit can be shadowed by a periodic orbit.
This proves that $p_n\in \cl(\per(X))$.
As the neighborhood $U$ is arbitrary,
we can assume that
$p_n\in \per(X)$ for every $n$.
Note that ${\cal O}_X(p_n)$ converges to
$\Lambda$ by Corollary \ref{cor1}.

As $\Lambda$ is transitive we have that
if $E^s\oplus E^X\oplus E^u$ denotes the
corresponding hyperbolic splitting,
then $\di(E^s)=s$ and  $\di(E^u)=u$ are constant
in $\Lambda$. Clearly neither
$s=0$ nor $u=0$ since $\Lambda$ is not
$\Omega$-isolated.
As ${\cal O}_X(p_n)$ converges to
$\Lambda$ both the local stable and unstable manifolds
of $p_n$ have dimension $s$, $u$ respectively.
Moreover, both invariant manifolds have
uniform size as well.
This implies that
$W^u_X(p_n)\cap W^s_X(p)\neq\emptyset$
and $W^s_X(p_n)\cap W^u_X(p)\neq\emptyset$
for $n$ large.
As $p_n\in Per(X)$ and $p\in \Lambda$,
we conclude by the Inclination Lemma
\cite{MP92} that
$p_n\in \cl(W^s_X(p)\cap W^u_X(p))$.
As $p\in \Lambda$, $W^{s,u}_X(p)\subset W^{s,u}_X(\Lambda)$.
So, $p_n\in \cl(W^s_X(\Lambda)\cap W^u_X(\Lambda))$.
As $\Lambda$ is saturated,
$W^s_X(\Lambda)\cap W^u_X(\Lambda)=\Lambda$
and hence
$p_n\in \cl(\Lambda)=\Lambda$.
But this is impossible since
$p_n\in \Omega(X)\setminus \Lambda$
by assumption.
This concludes the proof.
\qed

Denote by
$\mathcal{F}$ the collection of all 
isolated transitive neutral sets of $X$.

\begin{prop}
\label{p4}
A sub collection
$\mathcal{F}'$ of $\mathcal{F}$ is finite
if and only if
$\,\,
\cup_{\Lambda\in \mathcal{F}'} \Lambda
$
is closed.
\end{prop}

{\pf}
Obviously
$\cup_{\Lambda\in \mathcal{F}'} \Lambda$ is closed
if $\mathcal{F}'$
is finite.
Conversely, suppose that 
$\cup_{\Lambda\in \mathcal{F}'} \Lambda$ is closed.
If $\mathcal{F}'$ were infinite then,
it would exist sequence
$\Lambda_n\in \mathcal{F}'$ of (different) sets
accumulating some
$\Lambda\in \mathcal{F}'$.
By Corollary \ref{cor1}
we have
$\Lambda_n\subseteq U$ for some isolating
block $U$ of $\Lambda$ and $n$ large.
And then, we would have that
$\Lambda_n= \Lambda$ for $n$ large, a contradiction.
\qed

\section{Homoclinic classes}
\label{s-hc}

The main result of this section is
\begin{thm}
\label{h-neutral}
There is a residual subset ${\cal R}$ of ${\cal X}^1(M)$ such that
every homoclinic class of every vector field in $\cal R$ is neutral.
\end{thm}

\begin{clly}
\label{equivalence}
The following properties are
equivalent for $X\in {\cal R}$ 
and every compact invariant set $\Lambda$ of $X$ :
\begin{enumerate}
\item
$\Lambda$ is a transitive neutral set
with periodic orbits of $X$.
\item
$\Lambda$ is a homoclinic class of $X$.
\item
$\Lambda$ is a maximal transitive set with periodic
orbits of $X$.
\end{enumerate}
\end{clly}

{\pf}
That (2) implies (1) follows
from Theorem \ref{h-neutral}. 
That (1) implies (3) follows from Lemma \ref{l2}-2.
Let us prove that (3) implies (2).
If $\Lambda$ is as in (3) and $p\in \per(X)\cap \Lambda$, 
then $\Lambda\cap H_X(p)\neq \emptyset$.
By Theorem \ref{h-neutral} we can assume $H_X(p)$ is neutral,
and so it is maximal transitive (using $(1) \Rightarrow (3)$).
As both $\Lambda$ and $H_X(p)$ are maximal transitive
we conclude $\Lambda=H_X(p)$ and the proof follows.
\qed

\begin{clly}
\label{equival}
For $X\in {\cal R}$, a non singular
compact isolated set of $X$
is neutral and transitive
if and only if it is a homoclinic class.
\end{clly}

{\pf}
The converse follows from Theorem \ref{h-neutral}.
To prove the direct,
denote $\Lambda$ a transitive isolated
neutral set of a generic $C^1$ vector field $X$.
By Proposition \ref{p2} it follows that $\Lambda$
is also $\Omega$-isolated.
Since $\Lambda$ is transitive we have $\Lambda
\subset \Omega(X)$. Thus, by \cite{Pu67} it follows that
$\Lambda=\cl(\Lambda\cap \per(X))$,
and so, $\Lambda\cap \per(X)\neq\emptyset$.
Then the conclusion follows from the
previous corollary.
\qed

The proof of Theorem \ref{h-neutral}
follows immediately from the two lemmas below.

\begin{lemma}
\label{lm}
There exists a residual set $\mathcal{R}$ of $\mathcal{X}^1(M)$ such that,
for every $Y \in \mathcal{R}$ and $\sigma \in \cri(Y)$,
$\cl(W^u_X(\sigma))$ is Lyapunov stable for $X$ and
$\cl(W^s_X(\sigma))$ is Lyapunov stable for $-X$.
\end{lemma}

\begin{lemma}
\label{lh}
There exists a residual set $\mathcal{R}$ in $\mathcal{X}^1(M)$ such that
every $X \in \mathcal{R}$ satisfies
\[
H_X(p) = \cl(W^u_X(p)) \cap \cl(W^s_X(p))
\]
for all $p\in \per(X)$. 
\end{lemma}

Lemma \ref{lm} was proved in \cite[Theorem 6.1, p. 372]{MP01}
when $\sigma$ is a singularity and
the same proof works when $\sigma$ is a periodic orbit.
We shall give another proof of this lemma
in the Appendix for completeness.

Before the proof of Lemma \ref{lh}, let us introduce some notation.
Recall that $M$ is a closed $n$-manifold, $n\geq 3$.
We denote $2^M_c$ the space of all compact subsets
of $M$ endowed with the Hausdorff topology.
Recall that ${\cal KS}^1(M)\subset {\cal X}^1(M)$ denotes the set 
of Kupka--Smale $C^1$ vector fields on $M$. 

Given $X\in {\cal X}^1(M)$ and $p\in
\per(X)$ we denote $\Pi_X(p)$ the  period of $p$. 
We set $\Pi_X(p)=0$ if $p$ is a singularity of $X$. 

If $T>0$ we denote
$$
\cri_T(X)=\{p\in \cri(X):\Pi_X(p)<T\}.
$$

If $p\in \cri(X)$ is hyperbolic, then
there is a continuation $p(Y)$ of $p$
for $Y$ close enough to $X$ so that
$p(X)=p$.

Note that if $X\in {\cal KS}^1(M)$ and $T>0$,
then
$$
\cri_T(X)=\{p_1(X),\cdots,p_k(X)\}
$$
is a finite set.
Moreover,
$$
\cri_T(Y)=\{p_1(Y),\cdots,p_k(Y)\}
$$
for every $Y$ close enough to $X$.

Let ${\cal Y}$ be a metric space.
A set-valued map
$$
\Phi:{\cal Y}\to 2_c^M
$$
is {\em lower semi-continuous}
at $Y_0\in {\cal Y}$ if for every open set $U\subset M$ one has
$\Phi(Y_0)\cap U\neq\emptyset$ implies $\Phi(Y)\cap U\neq\emptyset$ for
every $Y$ close to $Y_0$.
Similarly, we say that $\Phi$ is {\em upper semi-continuous}
at $Y_1\in {\cal Y}$ if for every compact set $K\subset M$
one has $\Phi(Y_1)\cap K=\emptyset$ implies $\Phi(Y)\cap K=\emptyset$ 
for every $Y$ close to $Y_1$.
We say that $\Phi$ is {\em lower semi-continuous} if it is 
lower semi-continuous at every $Y_0\in {\cal Y}$.
A well known result \cite[Corollary 1, p. 71]{Ku68} asserts 
that if $\Phi:{\cal X}^1(M)\to 2_c^M$ is a 
lower semi-continuous map, then it is upper semi-continuous
at every $Y$ in a residual subset of
${\cal X}^1(M)$.

The lemma below is the flow version
of \cite[Theorem E, p. 5214]{WX00} 
(see also \cite{Ar99,Ha1,Ha97,Ha98}).

\begin{lemma}
\label{lwx}
Let $Y\in {\cal X}^1(M)$ and $x \notin \cri(Y)$.
For any $C^1$ neighborhood ${\cal U}$ of $Y$ there are
$\rho > 1$, $L >0$ and $\epsilon_0>0$ such that for any
$0 < \epsilon \leq \epsilon_0$ and any two points
$p,q\in M$ satisfying
\begin{description}
\item{(a)}
$p,q\notin B_\epsilon(Y_{[-L,0]}(x))$,
\item{(b)}
${\cal O}^+_Y(p)\cap B_{\epsilon/\rho}(x)\neq \emptyset$, and
\item{(c)}
${\cal O}^-_Y(q)\cap B_{\epsilon/\rho}(x)\neq \emptyset$,
\end{description}
there is $Z\in {\cal U}$
such that $Z=Y$ off $B_\epsilon(Y_{[-L,0]}(x))$
and that $q\in {\cal O}_Z^+(p)$.
\end{lemma}

{\flushleft{\bf Proof of Lemma \ref{lh}: }}
Given $X \in {\cal X}^1(M)$ we denote by $\per_T(X)$ the set of 
periodic orbits of $X$ with period $< T$.

We first prove a local version of Lemma \ref{lh}:
\begin{lemma}
\label{ll1}
If $X \in {\cal K}{\cal S}^1(M)$ and $T>0$ then there are a 
neighborhood ${\cal V}_{X,T} \ni X$ and a residual subset 
${\cal P}_{X,T}$ of ${\cal V}_{X,T}$ such that if 
$Y \in {\cal P}_{X,T}$ and $p\in \per_T(Y)$ then 
$H_Y(p)=\cl(W^u_Y(p))\cap \cl (W^s_Y(p))$.
\end{lemma}

\pf
There is a neighborhood ${\cal V}_{X,T}\ni X$ such that
$$
\per_T(Y)= \{\sigma_1(Y),\dots,\sigma_m(Y)\}\quad
\forall \quad Y\in {\cal V}_{X,T}.
$$
For each $1\le i \le m$, let 
$\Psi_i:{\cal V}_{X,T}\ni Y \mapsto H_Y(\sigma_i(Y))\in 2_c^M$.
Note that $\Psi_i, \forall \, i$, is lower semi-continuous by the 
persistence of transverse homoclinic orbits.
So, there is a residual subset ${\cal P}^i_{X,T}$ of ${\cal V}_{X,T}$ such
that $\Psi_i$ is upper semi-continuous in ${\cal P}^i_{X,T}$.
Set ${\cal P}_{X,T}={\cal K}{\cal S}^1(M) \cap(\cap {\cal P}^i_{X,T})
\cap {\cal R}$, where ${\cal R}$ is the residual set given in
Lemma \ref{lm}.
Then ${\cal P}_{X,T}$ is residual in ${\cal V}_{X,T}$.

Let us prove that ${\cal P}_{X,T}$ satisfies the conclusion of the lemma.
For this, let $\sigma\in \per_T(Y)$ for some $Y \in {\cal P}_{X,T}$.
Then $\sigma=\sigma_i(Y)$ for some $i$, and so $\Psi_i(Y)= H_Y(\sigma)$.

Suppose, by contradiction, that 
$H_Y(\sigma)\neq \cl (W^u_Y(\sigma))\cap \cl(W^s_Y(\sigma))$. 
Then there is 
$x\in \cl(W^u_Y(\sigma))\cap \cl(W^s_Y(\sigma))\setminus H_Y(\sigma)$.

We have either
\begin{description}
\item{(a)}
$x \notin \cri(Y)$ or
\item{(b)}
$x \in \cri(Y)$.
\end{description}

It is enough to prove the lemma in case (a).
Indeed, suppose that case (b) holds.
As $Y$ is Kupka--Smale we have that
${\cal O}_Y(x)$ is hyperbolic.
Clearly ${\cal O}_Y(x)$ is neither a sink or a source
and so $W^s_Y(x)\setminus {\cal O}_Y(x)\neq\emptyset$
and $W^u_Y(x)\setminus {\cal O}_Y(x)\neq
\emptyset$.
Note that $\cl(W^u_Y(\sigma))$ is Lyapunov stable
since $Y\in {\cal R}$. As $x\in \cl(W^u_Y(\sigma))$
we conclude that 
$W^u_Y(x)\subseteq \cl(W^u_Y(\sigma))$.
As $x \in \cl(W^s_Y(\sigma))$, there is
$x'\in \cl(W^s_Y(\sigma))\cap (W^u_Y(x)\setminus \mathcal{O}_Y(x))$
arbitrarily close to $x$
(for this use the Grobman--Hartman Theorem
as in the proof of Lemma \ref{lm}
in the Appendix).
Obviously $x'\notin \cri(Y)$.
If $x'\in H_Y(\sigma)$ we would have
that $x\in H_Y(\sigma)$ since $\alpha_Y(x')= \mathcal{O}_Y(x)$
contradicting  $x\notin H_Y(\sigma)$.
Henceforth 
$x'\in \cl(W^u_Y(\sigma))\cap \cl(W^s_Y(\sigma))\setminus H_Y(\sigma)$
and $x'\notin \cri(Y)$.
Then we conclude as in case (a)
replacing $x$ by $x'$.

Now we prove the lemma in case (a).

As $x \notin H_Y(\sigma)$, there is 
a compact neighborhood $K$ of $x$ such that
$
K\cap H_Y(\sigma)= \emptyset
$.
As $\Psi_i$ is upper semi-continuous at $Y$, 
there is a neighborhood ${\cal U}$ of $Y$ such that
\begin{equation}
\label{equino}
K\cap H_Z(\sigma(Z))= \emptyset,
\end{equation}
for all $Z\in {\cal U}$.
 
Let $\rho,L,\epsilon_0$ be the constants in Lemma \ref{lwx}
for $Y\in {\cal X}^1(M)$, $x$, and ${\cal U}$ as above.

As $x\notin \cri(Y)$,
$Y_{[-L,0]}(x)\cap {\cal O}_Y(\sigma)= \emptyset$.
Then, there is $0<\epsilon<\epsilon_0$ such that
${\cal O}_Y(\sigma)\cap B_\epsilon(Y_{[-L,0]}(x))=\emptyset$
and $B_\epsilon(x)\subseteq K$.

Choose an open set
$V$ containing ${\cal O}_Y(\sigma)$ such that 
$V\cap B_\epsilon(Y_{[-L,0]}(x))= \emptyset$.

As $x\in \cl(W^u_Y(\sigma))$, one can choose 
$p\in W^u_Y(\sigma)\setminus \{\sigma\}\cap V$
such that
$$
{\cal O}_Y^+(p)\cap B_{\epsilon/\rho}(x)\neq \emptyset.
$$

Similarly, as $x\in \cl(W^s_Y(\sigma))$, one can choose 
$q\in W^s_Y(\sigma)\setminus \{\sigma\}\cap V$
such that
$$
{\cal O}_Y^-(q)\cap B_{\epsilon/\rho}(x)\neq \emptyset.
$$

We can assume that ${\cal O}_Y^-(p) \subset V$ and
${\cal O}_Y^+(q) \subset V$.
Henceforth
\begin{equation}
\label{epi}
({\cal O}_Y^-(p)\cup {\cal O}_Y^+(q))
\cap B_\epsilon(Y_{[-L,0]}(x))= \emptyset.
\end{equation}

Observe that
$$
q\notin {\cal O}_Y^+(p)
$$
for, otherwise,
$p$ would be a homoclinic orbit of $Y$ passing through
$K$ contradicting (\ref{equino}).

By construction $\epsilon,p,q$ satisfy (b) and (c) of Lemma \ref{lwx}.

As $p,q\in V$ and $V\cap B_\epsilon(Y_{[-L,0]}(x))= \emptyset$
we have that that $\epsilon,p,q$ also
satisfy (a) of Lemma \ref{lwx}.

Then, by Lemma \ref{lwx},
there is $Z\in {\cal U}$ such that
$Z=Y$ off $B_\epsilon(Y_{[-L,0]}(x))$ and $q\in {\cal O}_Z^+(p)$.

Clearly $\sigma(Z)=\sigma$ and by (\ref{epi})
we have $p\in W^u_Z(\sigma)$ and $q\in W^s_Z(\sigma)$
since $Z=Y$ off $B_\epsilon(Y_{[-L,0]}(x))$.

Hence
${\cal O}= {\cal O}_Z(p)= {\cal O}_Z(q)$
is a homoclinic orbit of $\sigma$.

As $q\notin {\cal O}_Y^+(p)$, we have that
${\cal O}\cap B_\epsilon(x)\neq\emptyset$.

Perturbing $Z$ we can assume that ${\cal O}$ is transverse, 
i.e. ${\cal O}\subseteq H_Z(\sigma)$.

As ${\cal O}\cap B_\epsilon(x)\neq\emptyset$
and $B_\epsilon(x)\subset K$ we would obtain
$K\cap H_Z(\sigma(Z))\neq\emptyset$ contradicting
(\ref{equino}).

This finishes the proof.
\qed

{\flushleft{\bf Proof of Lemma \ref{lh}: }}
Fix $T>0$.
For any $X\in {\cal KS}^1(M)$
consider ${\cal V}_{X,T}$ and ${\cal P}_{X,T}$
as in Lemma \ref{ll1}.
Choose a sequence $X^n\in {\cal KS}^1(M)$ such that 
$\{X^n: n\in I\!\! N\}$ is dense in ${\cal X}^1(M)$
(recall that ${\cal X}^1(M)$ is a separable metric space).
Denote ${\cal V}_{n,T}={\cal V}_{X^n,T}$ 
and ${\cal P}_{n,T}={\cal P}_{X^n,T}$.

Define
$$
{\cal O}^T=\cup_{n}{\cal V}_{n,T}
\,\,\,\,\,\,\mbox{and}\,\,\,\,\,\,
{\cal P}^T=\cup_{n}{\cal P}_{n,T}.
$$
Clearly ${\cal O}^T$ is open and dense in ${\cal X}^1(M)$.

We claim that ${\cal P}^T$ is residual in ${\cal O}^T$.
Indeed, for any $n$ there is a sequence $D_{k,n,T}$, $k\in I\!\! N$, 
such that
$$
{\cal P}_{n,T}=\cap_k D_{k,n,T},
$$
and $D_{k,n,T}$ is open and dense in ${\cal V}_{n,T}$ for any $k$.
As
$$
{\cal P}^T=\cup_n{\cal P}_{n,T}=
\cup_n(\cap_kD_{k,n,T})=
\cap_k(\cup_nD_{k,n,T})
$$
and $\cup_nD_{k,n,T}$ is open and dense in $\cup_n{\cal V}_{n,T}=
{\cal O}^T$
we conclude
that ${\cal P}^T$ is residual in ${\cal O}^T$.
This proves the claim.

In particular, ${\cal P}^T$ is residual in
${\cal X}^1(M)$ for every $T$. 
Set ${\cal P}=\cap_{N\in I\!\! N}{\cal P}^N$.
It follows that ${\cal P}$ is
residual in ${\cal X}^1(M)$.
Choose $X\in {\cal P}$,
$p\in \per(X)$ and $N_0\in I\!\! N$ bigger
than $\Pi_X(p)+1$.
By definition $X\in {\cal P}^{N_0}$,
and so, $X\in {\cal P}_{X^n,N_0}$ for some $n$.
As $N_0>\Pi_X(p)$ we have $p\in \per_{N_0}(X)$.
Then $H_X(p)=\cl(W^u_X(p))\cap \cl(W^s_X(p))$ by Lemma \ref{ll1}
applied to $X^n$ and $T=N_0$.
This completes the proof of the lemma.
\qed

\section{Proof of Theorem \ref{thA} and Corollary \ref{c2}}
\label{s-pr}

{\flushleft{\bf Proof of Theorem \ref{thA}: }}
By Theorem \ref{h-neutral}
we have that homoclinic classes
for a residual subset of $C^1$ vector fields
on closed manifolds are neutral sets.
This leads us to apply the results
in Section \ref{s-ls}.
Thus,
Theorem \ref{thA}-(1) follows from Lemma
\ref{l1}-(2).
Theorem \ref{thA}-(2)
follows from Lemma \ref{l1}-(1).
Theorem \ref{thA}-(3) follows from
Corollary \ref{cor1}.
Theorem \ref{thA}-(4) follows from
Proposition \ref{p2}.
Similarly, Theorem \ref{thA}-(5) follows from
Proposition \ref{p2.5}.
Theorem \ref{thA}-(6) follows from
Proposition \ref{p3}

To prove Theorem \ref{thA}-(7)
we proceed as follows.
If $X$ has finitely many homoclinic classes,
then
the union of the homoclinic classes of
$X$ is obviously closed.
By \cite{Pu67}
it follows that $\Omega(X)$ is the union of the homoclinic
classes of $X$
(recall that $X$ is $C^1$ generic).
This implies that every homoclinic class of $X$ is 
$\Omega$-isolated, and so, they are
isolated by Theorem \ref{thA}-(4).
Conversely, suppose that
the union of the homoclinic classes
of $X$ is closed and
that every homoclinic class of $X$ is isolated.
Let $\mathcal{F}'$ be the collection of all
homoclinic classes of $X$.
As every homoclinic class of $X$ is isolated by
assumption one has $\mathcal{F}'\subset \mathcal{F}$
(recall the notation in Proposition \ref{p4}).
We have that $\cup_{\Lambda\in \mathcal{F}'}\Lambda$
is closed by hypothesis. Then,
Proposition \ref{p4} implies
that $\mathcal{F}'$ is finite and the proof follows.

{\flushleft{\bf Proof of Corollary \ref{c2}: }}
That (1) implies (2)
follows from Theorem \ref{thA}-(5).
To prove that (2) implies (1)
we proceed as follows.
If $H_X(p)$ is isolated, $H_X(p)$ is $\Omega$-isolated
by Theorem \ref{thA}-(4).
In particular, $H_X(p)$
is not in the closure of the sinks and sources
of $X$ unless it is either a sink or a source of
$X$ and we would be done. By \cite{Li83}, \cite{Man82},
as $M$ is 3-dimensional and
$X$ is nonsingular and generic,
one has that $H_X(p)$ is hyperbolic.
That (1) implies
(3) follows from the hyperbolic theory.
Indeed, if $H_X(p)$ is hyperbolic,
then $H_X(p)$ is isolated, transitive
and hyperbolic.
In other words, $H_X(p)$ is a basic set of $X$.
Then, the conclusion follows from the
structural stability
of basic sets \cite{PT93}.
That (3) implies (1) follows from \cite{Li83}, \cite{We95}
since $H_X(p)$ has no singularities
(recall $X$ has no singularities by hypothesis).
\qed

\section{Appendix}

Here we give a proof of Lemma \ref{lm} using Lemma \ref{lwx}.
The proof we gave in \cite[Theorem 6.1, p. 372]{MP01} uses a
different version of the $C^1$ Closing Lemma rather than Lemma \ref{lwx}.

As in the proof of Lemma \ref{lh}, 
Lemma \ref{lm} is a consequence of the following local lemma.

\begin{lemma}
\label{local-lm}
If $X\in {\cal KS}^1(M)$ and $T>0$, then there is a neighborhood 
${\cal U}_{X,T}$ of $X$ and a residual subset ${\cal R}_{X,T}$ of 
${\cal U}_{X,T}$ such that if $Y\in {\cal R}_{X,T}$ and 
$p\in \cri_T(Y)$, then $\cl(W^u_Y(p))$ is Lyapunov stable for $Y$
and $\cl(W^s_Y(p))$ is Lyapunov stable for $-Y$.
\end{lemma}

\pf
Recall that $\cri(Y)= \{p_1(Y),\dots,p_k(Y)\}$
for every $Y$ in some neighborhood ${\cal U}_{X,T}$ of $X$,
where $p_i(Y)$, $1\leq i \leq k$, is either a periodic orbit
or a singularity of $Y$.

For any $i\in \{1,\cdots, k\}$ we define
$
\Phi_i:{\cal U}_{X,T}\to 2_c^M
$
by
$$
\Phi_i(Y)=\cl(W^u_Y(p_i(Y)).
$$

By the continuous dependence of unstable manifolds
we have that $\Phi_i$ is a lower semi-continuous
map, and so, $\Phi_i$ is also upper semi-continuous
for every vector field in some residual subset 
${\cal R}_i$ of ${\cal U}_{X,T}$.
Set ${\cal R}_{X,T}= {\cal KS}^1(M)\cap (\cap_i {\cal R}_i)$.
Then ${\cal R}_{X,T}$ is residual in ${\cal U}_{X,T}$.
Let us prove that ${\cal R}_{X,T}$
satisfies the conclusion of the lemma.

Let $\sigma\in \cri_T(Y)$ for some
$Y\in {\cal R}_{X,T}$.
Then, $\sigma= p_i(Y)$ for some $i$, and so,
$\Phi_i(Y)= \cl(W^u_Y(\sigma))$.

Suppose by contradiction that
$\cl(W^u_Y(\sigma))$ is not Lyapunov stable for $Y$. 

Then, there are an open set
$U$ containing $\cl(W^u_Y(\sigma))$
and two sequences $x_n\to x\in \cl(W^u_Y(\sigma))$,
$t_n\geq 0$ 
such that
$$
Y_{t_n}(x_n)\notin U.
$$

As in the proof of Lemma \ref{lh} we have either
\begin{description}
\item{(a)}
$x \notin \cri(Y)$ or
\item{(b)}
$x \in \cri(Y)$.
\end{description}

Again it is enough to prove the lemma in case (a).
Indeed, suppose that case (b) holds.
As $Y$ is Kupka--Smale we have that ${\cal O}_Y(x)$ is hyperbolic.
Clearly ${\cal O}_Y(x)$ is neither a sink or a source
and so $W^s_Y(x)\setminus {\cal O}_Y(x)\neq\emptyset$
and $W^u_Y(x)\setminus {\cal O}_Y(x)\neq
\emptyset$.
Let $V \subset U$ be a small neighborhood of $x$ given by the 
Grobman--Hartman Theorem \cite{MP92} such that
$\partial(W^u_Y(x,V))= D^u_Y(x)$ is a fundamental domain for
$W^u_Y(x)$
(here $W^u_Y(x,V)$ denotes the connected component
of $V\cap W^u_Y(x)$ containing ${\cal O}_X(x)$).
Note that $D^u_Y(x)\subset W^u_Y(x)\setminus {\cal O}_Y(x)$. 
As $x_n\to x$, we can assume $x_n\in \interior(V)$ for all $n$.
As $Y_{t_n}(x_n)\notin U$, we have that $x_n\notin W_Y^s(\sigma)$. 
So, there is $s_n> 0$ such that $x_n'= Y_{s_n}(x_n)\in \partial V$ and
$Y_{s}(x_n)\in \interior(V)$ for $0\leq s < s_n$.
Since $Y_{t_n}(x_n)\notin U$ we have that
$Y_{t_n}(x_n)\notin \cl(V)$ for all $n$.
From this we conclude that $s_n< t_n$ for all $n$.
On the other hand, as $x_n\to x$,
passing to a subsequence if necessary,
we can assume that $x_n'\to x'$
for some $x' \in
D_Y^u(x) \subset W^u_Y(x)\setminus {\cal O}_Y(x)$.
Now we have the following claim.

\begin{claim}
$x'\in \cl(W_Y^u(\sigma))$.
\end{claim}

\pf
As $x\in \cl(W^u_Y(\sigma))$,
using the Connecting Lemma \cite{Ha97},
there is $Z$
$C^1$ near $Y$ such that 
$W_Z^u(\sigma(Z))\cap W^s_Z(x(Z))\neq \emptyset$.
In other words there is a saddle-connection
between $\sigma(Z)$ and $x(Z)$.
Breaking this saddle-connection
as in
the proof of Lemma 2.4 in \cite[p. 101]{MP92},
using the Inclination Lemma \cite{MP92},
we can find
$Z'$ $C^1$ close to $Z$
so that $W_{Z'}^u(\sigma(Z'))$ passes
close to $x'$.
This contradicts the upper-semicontinuity
of $\Phi_i$ at $Y$.
Thus, $x'\in \cl(W_Y^u(\sigma))$
and the Claim is proved.
\qed

As $x' \in D_Y^u(x)$ we have that $x'\notin \cri(Y)$.
As $x_n'\to x'$ and $Y_{t_n-s_n}(x_n') \notin U$
with $t_n-s_n >0$, we conclude as 
in case (a) replacing $x$ by $x'$,
$x_n$ by $x_n'$ and $t_n$ by $t_n-s_n$.

Now we prove the lemma in case (a).

As $\cl(W^u_Y(\sigma))\subseteq U$ and $\Phi_i$ is 
upper semi-continuous, there is a $C^1$ neighborhood 
${\cal U}\subset {\cal U}_{X,T}$ of $Y$ such that 
\begin{equation}
\label{equino1}
\cl(W^u_Z(\sigma(Z)))\subseteq U,
\end{equation}
for all $Z\in {\cal U}$.

Let $\rho, L, \epsilon_0$ as in Lemma \ref{lwx}
for $X=Y$, $x$, and ${\cal U}$ as above.

As $x\notin \cri(Y)$,
$Y_{[-L,0]}(x)\cap {\cal O}_Y(\sigma)= \emptyset$.

As $x \in \cl(W^u_Y(\sigma))$, $Y_{[-L,0]}(x)
\subset \cl(W^u_Y(\sigma))$ and so $Y_{[-L,0]}(x)\subset U$.
Then, there is $0<\epsilon\leq \epsilon_0$ such that
$B_\epsilon(Y_{[-L,0]}(x))\cap {\cal O}_Y(\sigma)=\emptyset$
and $B_\epsilon(Y_{[-L,0]}(x))\subseteq U$.

Choose an open set $V$ containing
${\cal O}_Y(\sigma)$,
$V\subset \cl(V)\subset U$,
such that 
$V\cap B_\epsilon(Y_{[-L,0]}(x))= \emptyset$.

For $n$ large, we have $x_n\in B_{\epsilon/\rho}(x)$ and we set $q= Y_{t_n}(x_n)\notin U$.

As $x_n\to x$ and $t_n>0$ we have ${\cal O}_Y^-(q)\cap B_{\epsilon/\rho}(x)
\neq \emptyset$.

As $x\in \cl(W^u_Y(\sigma))$, there is 
$p\in (W^u_Y(\sigma)\setminus \{\sigma\})\cap V$ such that
${\cal O}_Y^+(p)\cap B_{\epsilon/\rho}(x)\neq \emptyset$
and ${\cal O}_Y^-(p)\subseteq V$. 

By construction
$\epsilon,p,q$ satisfy (b) and (c) of Lemma \ref{lwx}.

As $V\cap B_\epsilon(Y_{[-L,0]}(x))=\emptyset$
and
$q\notin U$,
we have that $\epsilon,p,q$ also
satisfy (a) of Lemma \ref{lwx}.

Then, by Lemma \ref{lwx}, there is $Z\in {\cal U}$ such that
$Z= Y$ off $B_\epsilon(Y_{[-L,0]}(x))$ and
$q\in {\cal O}_Z^+(p)$.

As $V\cap B_\epsilon(Y_{[-L,0]}(x))= \emptyset$
and ${\cal O}_Y^-(p)\subset V$ we have
\begin{equation}
\label{equi1}
{\cal O}^-_Y(p)\cap
B_\epsilon(Y_{[-L,0]}(x))= \emptyset.
\end{equation}

Now, (\ref{equi1}), together with
$V\cap B_\epsilon(Y_{[-L,0]}(x))= \emptyset$
and $Z=Y$ off $B_\epsilon(Y_{[-L,0]}(x))$ imply that
$\sigma(Z)=\sigma$ and $p\in W_Z^u(\sigma)$.
As $ q\notin U$
and $q\in W^u_Z(\sigma)$
(recall $p\in W^u_Z(\sigma)$ and $q\in {\cal O}^+_Z(p)$),
we have a contradiction by (\ref{equino1}).
This finishes the proof.
\qed


\vspace{0.2cm}
\noindent C. M. Carballo\\
Departamento de Matem\'atica \\
PUC-Rio  \\
Rua Marqu\^es de S\~ao Vicente, 225\\
CEP 22453-900, 
Rio de Janeiro, R. J. , Brazil \\
{\em e-mail: carballo@mat.puc-rio.br}

\vspace{0.2cm}
\noindent C. A. Morales, M. J. Pacifico\\
Instituto de Matem\'atica \\
Universidade Federal do Rio de Janeiro \\
C. P. 68.530, CEP 21.945-970 \\
Rio de Janeiro, R. J. , Brazil \\
{\em e-mail: morales@impa.br}, \,
{\em pacifico@impa.br}

\end{document}